\documentclass[11pt]{article}
\usepackage{epsfig}
\usepackage{amsmath}
\usepackage{amsfonts}
\usepackage{amssymb}

\newtheorem{defn}{Definition}[section]

\newtheorem{thm}[defn]{Theorem}
\newtheorem{cor}[defn]{Corollary}

\newtheorem{prop}[defn]{Proposition}

\newcommand{\D}[1]{{\mathbb#1}}% Doubled -Blackboard bold - caps only

\newcommand{\RR}{{\D{R}}}

\begin{document}

\title{Basins of attraction for cascading maps}

\author{Erik Boczko\\
Department of Biomedical Informatics\\
Vanderbilt University\\
Nashville, TN 37232\\
615-936-6668\\
erik.m.boczko@vanderbilt.edu
\and Todd Young\\
Department of Mathematics\\
Ohio University\\
Athens, OH 45701\\
740-593-1285\\
young@math.ohiou.edu
}

%\thanks{}

%\keywords{}
%\subjclass{}

\maketitle

\begin{abstract}
We study a finite uni-directional array of ``cascading" or ``threshold coupled" 
chaotic maps. We describe some of the attractors for such systems and prove general
results about their basins of attraction. In particular, we show that
the basins of attraction have infinitely many path
components. We show that these components always accumulate 
at the corners of the domain of the system. For all threshold parameters above 
a certain value, we show that they accumulate at 
a Cantor set in the interior of the domain.
For certain ranges of the threshold, we prove that the system has many attractors.
\end{abstract}

\section{Cascading maps and applications}

Let $f(x) = 4x(1-x)$.  Let $c_1 \in (0,1)$ and call it a {\em threshold}. 
Formally, we define a {\em cascading}, or, {\em threshold coupled} system of $N$ maps as a
dynamical system on $[0,1]^N$ represented by a function 
$F_{c_1}:[0,1]^N \rightarrow [0,1]^N \times \RR$, 
that consists of iteration at each site by $f$, followed by ``cascading", i.e.,
$$
    F_{c_1} = \mathcal{C} \circ (f \times \cdots \times f),
$$
where $\mathcal{C}$ is a cascading operator.
Specifically,
$$
f \times \cdots \times f: ((x^j_1, \ldots, x^j_N)) \mapsto
               (y_1^j, \ldots, y_N^j) =  (f(x^j_1), \ldots, f(x^j_N))
$$
is simultaneous, independent iteration by $f$ at each site and 
\begin{equation}\label{cascade}
\begin{split}
     \mathcal{C} : &[0,1]^N \rightarrow [0,1]^N \times \RR  \\
                 : & (y^j_1, \ldots, y^j_N) \mapsto (x^{j+1}_1, \ldots, x^{j+1}_N, e^{j+1}).
\end{split}
\end{equation}
The map $F_{c_1}$ produces both the new state $(x^{j+1}_1, \ldots, x^{j+1}_N)$ 
as well as an {\em excess} $e^{j+1}$. The cascade map $\mathcal{C}$ acts on
$(y^j_1, \ldots, y^j_N)$ sequentially from left to right in the following fashion.
If $y^j_1 \le c_1$, then let $x^{j+1}_1 = y^j_1$ and define $e^j_1 = 0$.
If however, $y^j_1 \ge c_1$, then let $x^{j+1}_1 = c_1$ and define 
$e^j_1 = y^j_1 - c_1$. Next we consider the second site in the array in the same
way except we first add the value of $e_1^j$ to $y^j_2$. That is,
if $\hat{y}^j_2 = y^j_2 + e_1^j \le c_1$, then let $x^{j+1}_2 = \hat{y}^j_2$ 
and let $e^j_2 = 0$. Otherwise, let  $x^{j+1}_2 = c_1$ and define 
$e^j_2 = \hat{y}^j_1 - c_1$. The value of $e^j_2$ is cascaded forward to the third
site and this pattern continues to the end of the array. The final excess $e_N^j$
is the  $e^{j+1}$ that appears in (\ref{cascade}). It is 
``removed" from the system and the sequence $\{e^j\}$ is considered as the ``output". 

In applications, $c_1$ has been used as a tunable parameter and 
$\{e^j\}$ has been used as a time-series that characterizes the dynamics
of the system. Cascading maps have been considered for use in chaos-based computation 
\cite{SinDit99} and for classification of gene expression data \cite{Par02}. In these
potential applications initial conditions are classified based on the 
eventual output $\{e^j\}$, which depends on the attractor into which the orbit falls.
Thus it is important to have a clearer understanding of the attractors 
and their basins of attraction. An attractive aspect of cascading systems is that 
they can be applied using a single tunable parameter. 
Cascading dynamics of continuous time oscillators has also been considered \cite{Kap76}
and its use in communications has been explored \cite{KMS99}.

We will study attractors for $F$  and the associated basins, i.e. the set
of initial conditions which approach the attractors. In numerical studies
it appears that these basins are usually quite complicated. 
We prove that under some conditions, the basins of attraction have 
infinitely many components.
Our main result Theorem~\ref{int} shows that the components of the 
basins accumulate at a Cantor set in the interior of the domain. That is,
for each point in the Cantor set, every neighborhood of the point intersects
infinitely many components of the basins.
By component of a set, we mean connected component. In $\RR^N$ connectedness
is the same as path connectedness, i.e.\ any two points in a component can be connected
by a continuous curve. We also show the existence of at least two types of periodic attractors,
and that there are large ranges of threshold values
for which there many attractors.

For convenience we also define the related {\it threshold map} as
a map $f_{c_1}: [0,1] \rightarrow [0,1] \times \RR$ by
$$
f_{c_1}(x) = \left\{ \begin{array}{cl}
                      f(x), & \text{if } f(x) < c_1 \\
                      c_1, & \text{if } f(x) \ge c_1,
                   \end{array}
            \right.
$$
and we will denote the {\em excess} $f(x) - c_1$ for $f_{c_1}$ by $e(x)$.

%%%%%%%%%%%%%%%%%%%%%%%%%%%%%%%%%%%%%%%%%%%%%%%%%%%%%%%%%%%%%%%%%%%%%
%%%%%%%%%%%%%%%%%%%%%%%%%%%%%%%%%%%%%%%%%%%%%%%%%%%%%%%%%%%%%%%%%%%%%
%%%%%%%%%%%%%%%%%%%%%%%%%%%%%%%%%%%%%%%%%%%%%%%%%%%%%%%%%%%%%%%%%%%%%
%%%%%%%%%%%%%%%%%%%%%%%%%%%%%%%%%%%%%%%%%%%%%%%%%%%%%%%%%%%%%%%%%%%%%

\section{Dynamics of a single threshold map}

It is well known that $f = 4x(1-x)$ is topologically conjugate to the tent map.
(See for instance \cite[p.~114]{ASY97}.)
We will rely heavily on this fact in the proofs. Note that this method 
is not restrictive, but is quite general since any unimodal map $f$ with 
topological entropy $h_{top}(f)$ is semi-conjugate to a tent map with slope
$s = \exp(h_{top}(f))$ \cite{MilThu88}. In cases where $f$ is smooth and 
chaotic in the sense of having an absolutely continuous invariant measure, this
semi-conjugacy is in fact a conjugacy and smooth at almost all points. 
Also note that any orbit of $f$ that is bounded away from the critical 
point is uniformly hyperbolic (repelling). To see this, consider 
the metric provided by the conjugacy with the tent map. In that
metric the uniform expansion rate is 2. Finally, note that orbits of $f_{c_1}$ 
that stay below the threshold inherit these properties. For general maps $f$
under some restrictions, if an orbit is bounded away from a critical point
then it is hyperbolic (see Theorem III.5.1 in \cite{MelStr93}).

Let $C$ be the set of $x$ such that $f(x) \ge c_1$. One easily finds that
$$
  C = [c_0, 1-c_0] = \left[ \frac{1}{2} - \frac{1}{2}\sqrt{1 - c_1},  
                         \frac{1}{2} + \frac{1}{2}\sqrt{1 - c_1} \right].
$$
Thus $f_{c_1}(C) = \{c_1\}$ and $f^{-1}_{c_1}(c_1) = C$. 
Denote the forward orbit of $c_1$ by $c_2$, $c_3$, \ldots.

\begin{figure}
\centerline{\hbox{\psfig{figure=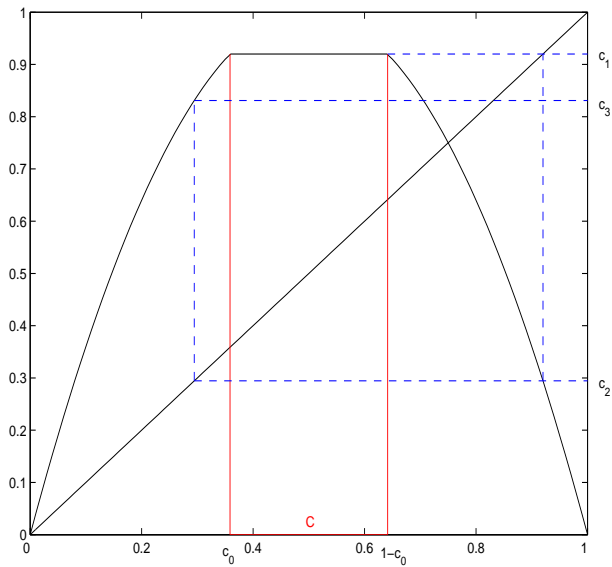,height=3in,width=3.2in}}}
\caption{Three iterations of the threshold map $f_{0.92}$. The interval
$C = [c_0,1-c_0]$ is the set on which $f_{c_1}(x) = c_1$. Almost all 
orbits eventually enter $C$. The interval $A = [c_2,c_1]$ absorbs all 
points except $0$ and $1$.}
\label{fig1}
\end{figure}

For $c_1 < 3/4$, it is easily seen that $f_{c_1}$ has a single attracting 
periodic point. We will restrict our attention to $c_1 > 3/4$.
The following proposition is elementary.
\begin{prop}
For $c_1 > 3/4$, $f_{c_1}$ has an absorbing set $A = [c_2, c_1]$ and
all points, except 0 and 1, are eventually mapped into $A$. 
\end{prop}

Given a periodic orbit, the set of all points that have backward orbits
which are asymptotic to the orbit of $x$ we call the {\em unstable manifold}
of the periodic orbit.
\begin{prop}\label{Wu34}
For $c_1 > 3/4$, the unstable manifold of  the fixed point at $3/4$ is $A$.
Given any point $a$ in $A$, there is a sequence of points
$\{a_{-j}\}_{j=1}^\infty$ that converges to $3/4$ and such that $f^j_{c_1}(a_{-j}) = a$.
Further, $a_{-j}$ can be chosen so that $\{a_{-i}\}_{i=1}^j$ does not intersect $C$.
\end{prop}
{\em Proof:} For $c_1 > 3/4$ it is easily seen that $f([1-c_0,c_1]) = A$ and that
$[1-c_0,c_1]$ is in the unstable manifold of $3/4$. Since $a$ has a preimage in $[1-c_0,c_1]$ 
and since $f$ is invertible on 
$[1-c_0,c_1]$ all preimages of $a$ exist in this interval. 
Since $(1-c_0,c_1]$ is mapped onto all of $A$ except for $c_1$
iterates of $a_{-j}$ may avoid $C$. \hfill $\Box$

\begin{prop}\label{Wuzero}
For $c_1 >3/4$, the unstable manifold of the fixed point at $0$ is $[0,c_1]$.
Given any point $a$ in $A$, there is an integer $j_0$ and a sequence of points
$\{a_{-j}\}_{j=j_0}^\infty$ that converges to $0$ such that $f^j_{c_1}(a_{-j}) = a$.
\end{prop}
{\em Proof:} Note that $[0,c_0]$ is in the unstable manifold of $0$ and that
its image is $[0,c_1]$. The proof is then as for the previous proposition. \hfill $\Box$

Outside of $A$ the only invariant set is the fixed point at $0$.

Inside $A$ we may have two types of invariant sets, a super-stable or semi-stable
periodic orbit, and a hyperbolic repelling set $\Lambda$ in $A \setminus C$. 
By super-stable we mean that the multiplier of the periodic 
orbit is $0$, i.e.\ the orbit passes through the interior of $C$.
A semi-stable orbit occurs at threshold values for which a periodic orbit 
intersects the boundary of $C$.  Because of the hyperbolicity 
outside $C$, there are only two possibilities for the forward orbit of $c_1$; it 
lands in $C$ or it lands on $\Lambda$.

It follows from properties of $f$ that any hyperbolic repelling 
set $\Lambda$ must have measure zero  and any super-stable orbit
must attract Lebesgue almost all other points. Both of these conclusions follow from the 
fact that the forward orbit of Lebesgue almost every point must intersect $C$. This follows
easily for the case $f=4x(1-x)$ by the conjugation with the tent map. For more general maps
see Theorem 2.2 in \cite{HomYou02}.

The dynamics of the threshold map is very simple for 
$\frac{3}{4} < c_1 < \frac{5 + \sqrt{5}}{8} \approx 0.904508497$ as described by the next result.
\begin{prop}\label{trivial}
For 
$$
   \frac{3}{4} < c_1 < \frac{5 + \sqrt{5}}{8}
$$
$f_{c_1}$ has a super-stable period 2 orbit and $\Lambda = \{3/4\}$.
\end{prop}
{\it Proof:} One can easily solve $f(c_1) \in C$ to obtain the limits
of the above inequality. For such threshold values $c_2 > c_0$, so
$A \setminus C = [1-c_0,c_1]$. All points in this interval except $\{3/4\}$
eventually leave the interval and must then enter $C$. \hfill $\Box$ 

\begin{figure}
\centerline{\hbox{\psfig{figure=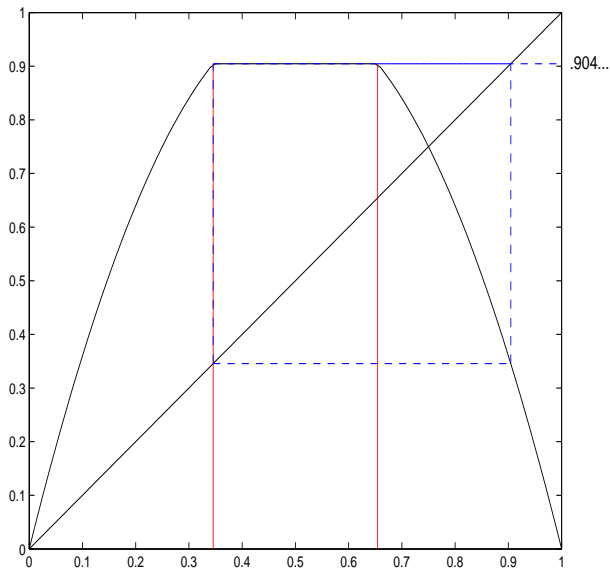,height=3in,width=3.2in}}}
\caption{For $c_1 = \frac{5+\sqrt{5}}{8} \approx .904508497$, $f_{c_1}$  maps
$c_1$ onto $c_0$. For all $c_1$ in  $(3/4,\frac{5+\sqrt{5}}{8})$, 
the threshold value $c_1$ is a super-stable period 2 point.}
\label{fig904}
\end{figure}

From numerical studies such as in \cite{Sin94} it is apparent that the bifurcations
become complex after $c_1 = \frac{5 + \sqrt{5}}{8}$. These bifurcations accumulate
at the value $c_1 = \xi_2 \equiv \frac{2 + \sqrt{3}}{4} \approx 0.933012702$, which is the location of
the second ``star" in the bifurcation diagram as described in \cite{Sin94}. In fact we 
have:
\begin{prop}\label{xi2bif}
There exists a sequence of open intervals $\{J_k\}_{k = 2}^\infty$ that accumulates at 
$\xi_2$, such that for $c_1 \in J_k$, the map has
a super-stable periodic orbit of period $k$. For $k$ even, $J_k$ is to the left of 
$\xi_2$ and $J_k$ is to the right of $\xi_2$ when $k$ odd.
\end{prop}
{\it Proof:}
The parameter value $\xi_2$ is a solution of 
$$
         f^2(\xi_2) = \frac{3}{4}.
$$
That is, the threshold value is mapped exactly onto the unstable fixed point at 3/4.
Note that there exist a sequence of inverse images of $C$ that accumulate on both sides of $3/4$. 
These inverse images vary smoothly with respect to $c_1$. Further note that as $c_1$ passes from 
$\frac{5 + \sqrt{5}}{8}$ to $1$, $f^2(c_1)$ passes from $1-c_0$ to $1$, Thus as $c_1$ is varied,
$f^2(c_1)$ must pass through all of the inverse images of $C$ surrounding $3/4$. 
As it does the super-stable orbits occur. \hfill $\Box$

\begin{figure}
\centerline{\hbox{\psfig{figure=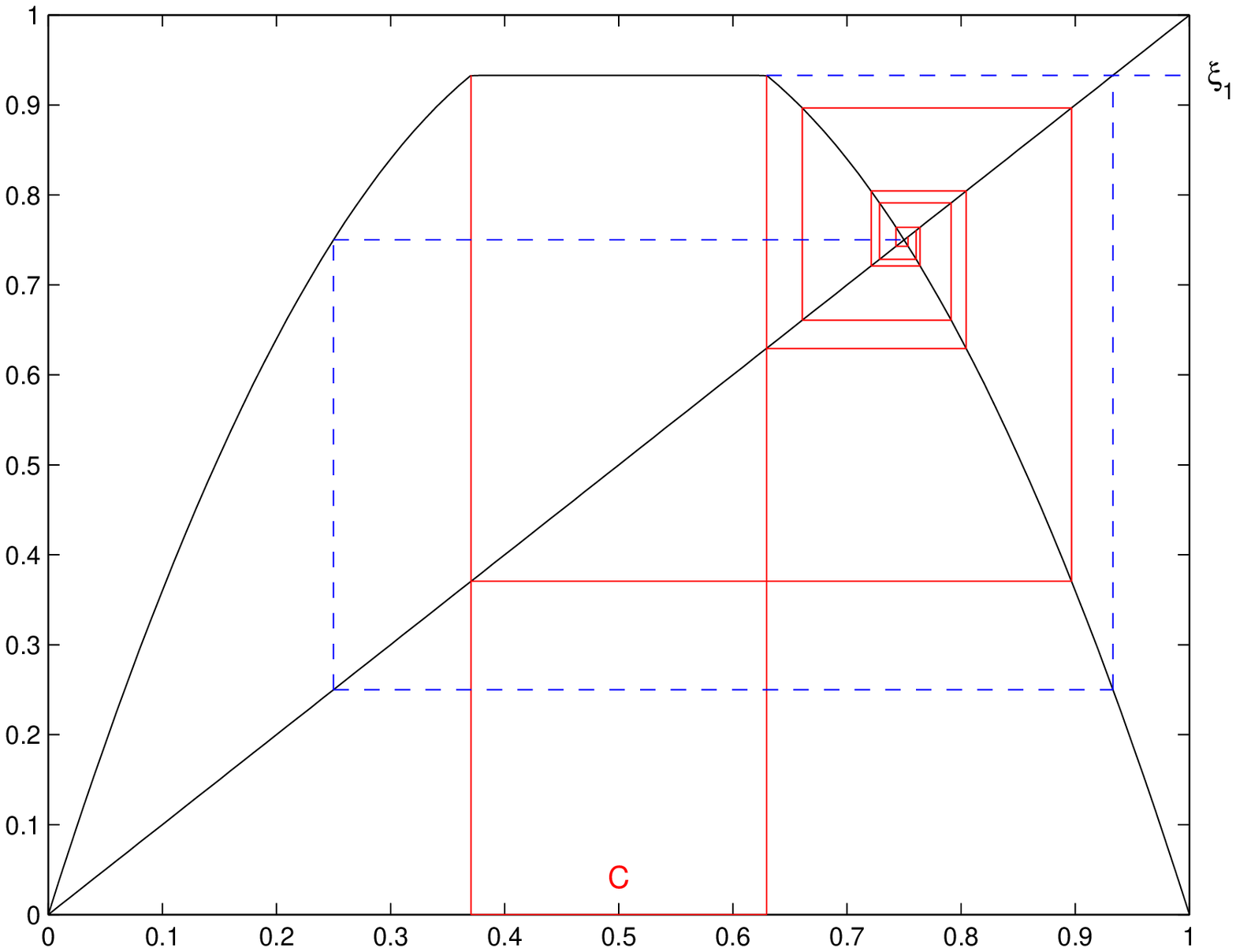,height=3in,width=3.2in}}}
\caption{The threshold map $f_{\xi_2}$ maps $C$ onto the unstable fixed point at $3/4$ in three iterations.
Also, inverse images of $C$ accumulate at $3/4$.}
\label{figxi2}
\end{figure}

Sinha \cite{Sin94} showed that there also exist threshold values $\xi_s$, $s>2$, approaching 1 for which
$$
   f^s_{\xi_s}(\xi_s) = \frac{3}{4}
$$
and $f^i_ {\xi_s}(\xi_s) < \frac{1}{2}$, for all $1 \le i \le s-1$. Further, he showed that
these values satisfy
$$
 \lim_{s \rightarrow \infty} \frac{\xi_{s+1} - \xi_s}{\xi_s - \xi_{s-1}} = \frac{1}{4}.
$$
These parameter values correspond to ``stars" in the bifurcation diagram. They can be understood
as highly degenerate homoclinic bifurcations. For each such parameter, the repelling fixed point at
3/4 has an interval of homoclinic orbits. Each of these threshold values has bifurcation structures
similar to those in Proposition~\ref{xi2bif}

Let $\{b_{2k}\}$ be the right endpoints of the even-index intervals $\{J_{2k}\}$.
\begin{prop}
For each $b_{2k}$, the map $f_{b_{2k}}$ has a semi-stable period $2k$ orbit. This orbit
persists and is repelling for all $c_1 > b_{2k}$.
\end{prop}  
{\it Proof:} For $b_{2k}$, the right endpoint $1-c_0$ of $C$ is $2k$ periodic. This periodic orbit
is a periodic orbit of $f$ that intersects $C$ only at $1-c_0$. For higher threshold values
$C$ is strictly smaller, so it does not intersect the periodic orbit. Thus the periodic
orbit persists for larger $c_1$. It is hyperbolic repelling since it does not intersect $C$. 
\hfill $\Box$

Thus  the hyperbolic attractor $\Lambda$ becomes increasingly complex as $c_1$ moves 
from   $\frac{5 + \sqrt{5}}{8}$ to $\xi_2 \equiv \frac{2 + \sqrt{3}}{4}$. At $\xi_2$ it 
becomes truly complex and the underlying dynamics become chaotic.

\begin{thm}\label{nontrivial}
For $c_1 > \xi_2$, $f_{c_1}$ has positive topological entropy  
and $\Lambda$ is uncountable. Further, $f$ restricted to $\Lambda$
is transitive, i.e. there is a heteroclinic orbit between any two 
periodic points.
\end{thm}
{\it Proof:}
For $c_1 < \xi_2$, $f_{c_1}((c_2,c_0))$ contains $3/4$, and at the same time
$(c_2,c_0) \subset A = W^u(3/4)$.  In particular, there exists a sequence 
of intervals $\{J_{-i}\}$, $i = 1,2,\ldots$, that are in the inverse image of
$J_0 = (c_2,c_0)$ that converge to $3/4$ as $i \rightarrow \infty$. Thus, 
for all $i \le - n_0$, we have that $J_{-i} \subset f_{c_1}(c_2,c_0)$. The 
collection $\{J_{-i}\}_{0}^{\infty}$ forms a Markov partition (see \cite{KatHas95}). 
If we consider a finite number of these intervals, $\{J_{-i}\}_{0}^{n}$, with $n \ge n_0+1$
then the transition matrix for this partition has zeros everywhere except for the
ones on super-diagonal and ones in the last $n - n_0 +1$ entries in the first column. 
$$
    M_n = \left( \begin{array}{ccccccc}
               0 & 1 & 0 & 0   & 0   & \cdots & 0 \\
               0 & 0 & 1 & 0   & 0   & \cdots & 0 \\
          \vdots & 0 & 0 & \ddots& 0  & \cdots & 0 \\ 
               0 & \vdots  & \ddots &    \ddots     & \ddots & \ddots & \vdots \\
               1 & \vdots &   \cdots  & \cdots  &    0  & 1 &  0 \\
          \vdots &  \vdots      & \cdots    &  \cdots &    0 & 0 &  1 \\
               1 & 0 & \cdots   &  \cdots  &  0 & 0 &  0
               \end{array} \right).
$$
These matrices have an eigenvalue greater than $1$, so we have
that there is an invariant set $\Lambda_n$ in $\{J_{-i}\}_{0}^n$ on which
$f_{c_1}$ is semi-conjugate to the Markov shift defined by $M_n$ (see Theorem 15.1.5 in 
\cite{KatHas95}). Hyperbolicity implies the semi-conjugacy is a conjugacy.
In particular $\Lambda_n$ is uncountable and $f_{c_1}$ restricted to it has
positive topological entropy and is transitive. Since $\Lambda_n$ is invariant, it is contained in
$\Lambda$. Further, $\Lambda = \left( \bigcup_{n = n_0}^\infty \Lambda_n \right) \bigcup \{3/4\}$ 
since a point in $\Lambda$ other than $3/4$ must enter $J_0$, every point in 
$\Lambda$ must be contained $J_0$ or one of the 
inverse images of $J_0$.   Since $f$ is transitive on each $\Lambda_n$, it is transitive
on $\Lambda$ which is the closure of the union of the set $\Lambda_n$.   \hfill $\Box$

The following is an easy consequence of Proposition~\ref{Wu34}
\begin{prop}\label{WuLat}
For $c_1 > \xi_2$, the unstable manifold of any periodic point in 
$\Lambda$ is $A$.
\end{prop}

Finally, we show that a super-stable orbit almost always occurs.
\begin{prop}\label{superstable}
For almost all $c_1$,  $f_{c_1}$ has a super-stable orbit.
\end{prop}
{\it Proof:}
Using the conjugacy with the tent map, one easily sees that the 
inverse images of $C$ form a full-measure subset of $[0,1]$. Its 
complement, which is the invariant set $\Lambda$, must have measure
zero.
Let $\mathcal{C}$ denote the collection of all intervals $J$ such 
that $f^k_{c_1}(J) = C$ and denote by $k(J)$, the minimal such integer $k$.
If $f_{c_1}$ does not have a super-stable orbit, then the orbit of $c_1$ must
be disjoint from $C$, so its forward orbit must be in the invariant set $\Lambda$, i.e.\
it must land on the hyperbolic repelling set.
This set has measure zero for all $c_1$ and points in the hyperbolic set
do not depend on $c_1$, since they are inherited from $f$. On the other hand, the 
points $c_k$ have a non-zero derivative with respect to $c_1$ as long as the 
orbit stays out of $C$. Thus for any $k$ the set of $c_1$ for which $c_k$ intersects the hyperbolic
set has measure zero. Taking the union of these sets for all $k$ we obtain a measure zero set.
The compliment of this measure zero set is the set of $c_1$ which never land on the 
hyperbolic set and thus must have a super-stable orbit. 
\hfill $\Box$

\begin{table}
\centerline{
\begin{tabular}{|c|c|c|} \hline
   $c_1$             & Attractor & $\Lambda$ \\ \hline
.75 -- .9045...      & period 2  & $\{3/4\}$ \\
.9045... -- .9330... & various   & expanding \\
.9330... -- 1.0      & various   & uncountable  \\ \hline
\end{tabular}
           }
\label{table1}
\caption{Summary of dynamics for a single threshold map.}
\end{table}

%%%%%%%%%%%%%%%%%%%%%%%%%%%%%%%%%%%%%%%%%%%%%%%%%%%%%%%%%%%%%%%%%%%%%
%%%%%%%%%%%%%%%%%%%%%%%%%%%%%%%%%%%%%%%%%%%%%%%%%%%%%%%%%%%%%%%%%%%%%
%%%%%%%%%%%%%%%%%%%%%%%%%%%%%%%%%%%%%%%%%%%%%%%%%%%%%%%%%%%%%%%%%%%%%
%%%%%%%%%%%%%%%%%%%%%%%%%%%%%%%%%%%%%%%%%%%%%%%%%%%%%%%%%%%%%%%%%%%%%

\section{Consequences for cascading}

Let $F_{c_1}$ be a sequence of $N$ threshold maps with cascading.
If one were to consider $N$ threshold maps without cascading, then
the results of the last section would give a very detailed picture
of the dynamics. Almost all initial conditions would fall into
superstable periodic orbits at each site. These periodic orbits
would exist at any  phases relative to one another. Thus, there
would be $p^{N-1}$ distinguishable attractors. Basins of the 
attractors would be cross-products of the inverse images of $C$
at each site, i.e.\ they would consist of rectangular boxes.
One might then consider cascading as a (possibly not small) 
perturbation. We show in the next section that many features
of the unperturbed system persist.

\subsection{Absorbing sets and invariant sets}

\begin{prop}
For $c_1 > 3/4$, $F_{c_1}$ has an absorbing set $A^N = A \times \cdots \times A$
that attracts almost all orbits.
\end{prop}
{\it Proof:}  For the threshold map $f_{c_1}$, $A$ is an absorbing set that 
attracts almost all orbits, and thus almost any initial condition in the first 
site will reach $A$. For sites after the first, consider a value $x_i^0$ at 
the $i$-th site. If $x_i^0 < c_0$, then the dynamics, including the effects 
of cascading only increase $x_i$, but cannot increase it past $c_1$. 
Since there are no fixed points in $[0,c_0]$ it will eventually
enter $[c_0,c-1]$ and at least on the next iteration will be in $A$. 
If $x_i^0 > c_1$, then $x_i^1 < c_0$. 

Once the orbit reaches $A$ in any given map, cascading cannot remove it
from $A$, since cascading can only increase the value at an individual map, but
cannot increase it past $c_1$ which is the upper bound of $A$.\hfill $\Box$

\begin{prop}\label{LambdaN}
For $c_1 > 3/4$,  the maximal hyperbolic repelling set of $F_{c_1}$ in $A^N$ is
$\Lambda^N = \Lambda \times \cdots \times \Lambda$.
\end{prop}
{\it Proof:} Suppose that $x_i^j \in \Lambda$ for each $i$, Then
the dynamics proceed at each site without cascading. On the other hand
if $(x_1, \ldots, x_N)$ is in the maximal repelling set, then the forward orbit of $x_i$
under $f_{c_1}$ cannot leave $A \setminus C$. For the first site which is not subject 
to the effects of cascading, this implies that $x_1$ is in $\Lambda$. This implies
that the first site will not have any excess, so the  second site is not subject
to cascading. Induction proves the result. \hfill $\Box$

The next results follow from Proposition~\ref{LambdaN} and Propositions~\ref{trivial} 
and \ref{nontrivial}.
\begin{cor}
For $3/4 < c_1 < \frac{5 + \sqrt{5}}{8} \approx .9045...$ the invariant set $\Lambda^N = \{3/4\}^N.$
\end{cor}
\begin{cor}
For $c_1 > \xi_2 \approx .9336...$ the invariant set $\Lambda^N$ is uncountable.
\end{cor}

\subsection{Structure of the basins}

\begin{thm}\label{bound}
The basins of all attractors accumulate at the corners of $[0,1]^N$. If a basin 
has a component contained in the interior of $[0,c_1]$ then every neighborhood
of every corner of $[0,1]^N$ contains infinitely
many components of the basin.
\end{thm}
{\it Proof:} Let $(x_1^0, \ldots, x_N^0)$ be a point in the basin of an attractor.
By Proposition~\ref{Wuzero}, we can find a point $(x^{-j}_1, \ldots, x^{-j}_N)$
near $(0, \ldots,0)$ such that $f_{c_1}^j(x_i^{-j}) = x_i^0$ in such a way that the orbit is disjoint
from $C$. Thus $(x^{-j}_1, \ldots, x^{-j}_N)$ is mapped onto $(x_1^0, \ldots, x_N^0)$ by $j$ $j$ 
iterations of $F$. This shows that the basins accumulate at $(0, \ldots, 0)$. 
If a basin has a component in $[0,c_1]$, it must have a sequence of inverse images
converging to $(0,\ldots, 0)$. Each of those inverse images must be components,
otherwise the original component could be extended along a continuous path.
Since $f(1 -x ) = f(x)$
this also implies the results at all the corners of $[0,1]^N$. \hfill $\Box$

\begin{prop}\label{WuLatN}
For $c_1 > 3/4$ the unstable manifold of the fixed point $\{3/4\}^N$ is $A^N$.
\end{prop}
{\it Proof:} 
Let $a = (a_1, a_2, \ldots, a_N)$ a point in $A^N$.
It follows from Proposition~\ref{WuLat} that given any $j$ 
there are points $a^{-j}_i$ in $[1-c_0,c_1]$
are mapped onto $a_i$ by $j$ iterations of $f_{c_1}$.  Further, 
the first $j-1$ iterates of $a^{-j}_i$ do not intersect $C$, thus there
is no cascading and the dynamics proceeds independently at each site.
\hfill $\Box$

Let $\tilde{\Lambda}$ be  the transitive component of $\Lambda$
containing $\{3/4\}^N$. That is a point $x \in \Lambda$ is in
$\tilde{\Lambda}$ if and only any neighborhood of $x$ contains points
that are eventually mapped onto $\{3/4\}^N$.

Our main result is the following.
\begin{thm}\label{int}
The basins of all attractors accumulate
at each point in $\tilde{\Lambda}^N$. If a basin
has a component contained in $A^N$, then each neighborhood of
each point in $\tilde{\Lambda}^N$ contains infinitely many components.
For $c_1 > \xi_2$, the above statements are true for $\Lambda$. 
\end{thm}
{\it Proof:} 
It follows from Proposition~\ref{WuLatN} that given any point $x$
in $\tilde{\Lambda}^N$ and any point $a$ in $A^N$, 
there is a point arbitrarily close to $y$ that is mapped into $a$.
It then follows from the fact that 
$F_{c_1}$ has a one to one inverse on $[1 - c_0,1]^N$ that components accumulate
at $\{3/4\}^N$. The rest of the conclusions follow from Propositions~\ref{nontrivial}
and~\ref{WuLatN}. 
\hfill $\Box$

\subsection{Examples of attractors}

\begin{prop}\label{inphase}
Suppose that for $c_1$, the threshold map, $f_{c_1}$, has a super-stable periodic orbit of period $p$.
Then $F_{c_1}$ has a super-stable periodic orbit of period $p$ corresponding to in-phase
synchronization of all $N$ maps. The basin of attraction of this periodic orbit includes
$C^N$ and all its preimages.
\end{prop}
{\it Proof:} 
Suppose that the initial state $(x_1^0, x_2^0, \ldots, x_N^0)$ is
contained in $C^N$. Then it is clear that one iteration of the system yields
$(x_1^1, \ldots, x_N^1) = (c_1, \ldots, c_1)$ and an excess which is the sum
of the excesses of each map. Since $p$ is minimal $(x_1^j, x_2^j, \ldots, x_N^j) = (c_j, \ldots, c_j)$
for all $1 \le j < p$ and $(x_1^p, \ldots, x_N^p) = (c_1, \ldots, c_1)$. The excess
on step $p$ will then be $N$ times the excess of each of the individual maps.
Clearly the basin of attraction of the synchronized periodic orbit includes the preimages of $C^N$.
\hfill $\Box$

The next result follows immediately from Propositions~\ref{superstable}~and~\ref{inphase}.
\begin{cor}
For almost all $c_1$, $F_{c_1}$ has a super-stable periodic orbit. 
\end{cor}

Next we consider another type of super-stable periodic orbit.
We call a periodic orbit a {\em ripple} if
$x_i^j = x_{i+1}^{j+1}$, for any $1 \le i \le N -1$.

\begin{prop}
Suppose that $c_1 > 15/16$ and that $f_{c_1}$ has a super-stable period $p$ orbit.
Then $F_{c_1}$ has a super-stable, ripple periodic orbit. This orbit has a basin with 
a non-empty interior.
\end{prop}
{\it Proof:}
Suppose the initial state $(x_1^0, \ldots, x_N^0)$ is such that $x_1^0 \in C$,
and 
$$
        x_i^{i-1} < x_{i-1}^{i-2}
$$
for each $2 \le i \le N$. Clearly we can find an open set of such initial conditions by
taking $x^0_i$ in one of the $i$-th inverse images of $(c_0,c_{p-1})$.
On the first step the excess from $f(x_1^0)$ is added to $f(x_2^0)$. The condition
$c_1 >15/16$ implies that $|f'(x)| < 1$ for all $x \in C$. This implies that
$x_i^{i-1}$ is in $C$, and thus $x_i^i = c_1$. \hfill $\Box$

For $N = 2$ and period $p = 2$ a ripple orbit might also be called an
anti-phase periodic orbit. We can give a necessary and sufficient condition for the existence of
such orbits.
\begin{prop}\label{antiphase}
Suppose that $f_{c_1}$ has a super-stable period $2$ orbit and $N = 2$.
Then $F_{c_1}$ has a super-stable anti-phase orbit if and only if
\begin{equation}\label{nokick}
c_2 + e(c_2) \le 1-c_0.
\end{equation}
The basin of this periodic orbit has a non-empty interior.
\end{prop}
{\it Proof:} Consider the initial condition $(x_1^0,x_2^0) = (c_2,c_1)$. 
Since $f(c_1) = c_2$ and $f(c_2) = c_1 + e(c_2)$, we have that
$(x_1^1,x_2^1) = (c_1,c_2+e(c_2))$. It then follows that
$(x_1^2,x_2^2) = (c_2,c_1)$. Thus the initial condition is part 
of a super-stable anti-phase orbit. It is clear that all initial 
conditions in a small open neighborhood of $(c_2,c_1)$ will be mapped onto 
$(c_2,c_1)$ by $F^2_{c_2}$. 

Now suppose that the condition (\ref{nokick}) is
not satisfied. Any anti-phase orbit must contain the point $(x_1^j,x_2^j) = (c_2,c_1)$.
However, if we apply the cascading map to $(c_2,c_1)$ and $c_2 + e(c_2) > 1-c_0$
then $x_1^{j+1} \neq c_1$, since $c_2 + e(c_2) \not\in C$.
\hfill $\Box$.

For the only period two window  $3/4 < c_1 < \frac{5 + \sqrt{5}}{8}$,
the equation in Proposition~\ref{antiphase} is equivalent to:
\begin{equation}
19 c_1 - 84 c_1^2 + 128 c_1^3 - 64 c_1^4 < \frac{1 + \sqrt{1-c_1}}{2}.
\end{equation}
This can be solved symbolically. Using a numerical approximation the 
anti-phase orbit exists for 
\begin{equation}\label{twophase}
    0.83627234814318 \ldots < c_1 < \frac{5 + \sqrt{5}}{8} \approx .904508497.
\end{equation}

Sinha and Ditto \cite{SinDit99} considered the parameter value $c_1 \approx .84$
where both in-phase and anti-phase periodic orbit exist.

Finally, consider again the case of $N>2$. 
\begin{thm}\label{multiple}
Suppose that $N>2$ and $c_1$ is as in (\ref{twophase}) or $c_1 > 15/16$
and $f_{c_1}$ has a period $p$ super-stable orbit..
Then there exist at least $2^{N-1}$ or $p^{N-1}$ attractors, respectively.
\end{thm}
{\it Proof:} Since the dynamics at the first two sites is independent of the
rest of the sites, both in-phase and anti-phase orbits exist in the
first two sites. It follows from arguments above that in-phase and
ripple synchronization of the third site with the second site exist. 
The rest are obtained by induction. \hfill $\Box$.

%%%%%%%%%%%%%%%%%%%%%%%%%%%%%%%%%%%%%%%%%%%%%%%%%%%%%%%%%%%%%%%%%%%%%
%%%%%%%%%%%%%%%%%%%%%%%%%%%%%%%%%%%%%%%%%%%%%%%%%%%%%%%%%%%%%%%%%%%%%
%%%%%%%%%%%%%%%%%%%%%%%%%%%%%%%%%%%%%%%%%%%%%%%%%%%%%%%%%%%%%%%%%%%%%
%%%%%%%%%%%%%%%%%%%%%%%%%%%%%%%%%%%%%%%%%%%%%%%%%%%%%%%%%%%%%%%%%%%%%

\section{Observations and conclusions}

\begin{figure}
\centerline{\hbox{\psfig{figure=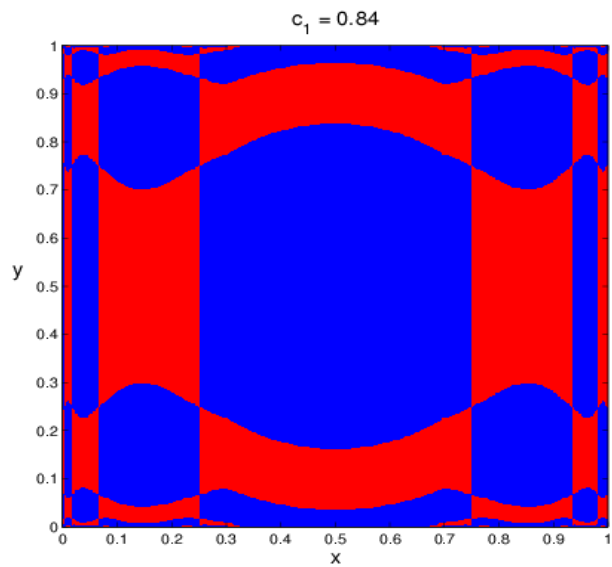,height=3in,width=3.2in}}}
\caption{Basins of attraction for $N=2$ and $c_1 = 0.84$. Here $c_1$
is in the range where both in-phase and anti-phase period 2 attractors 
exist.}
\label{figc84}
\end{figure}

\begin{figure}
\centerline{\hbox{\psfig{figure=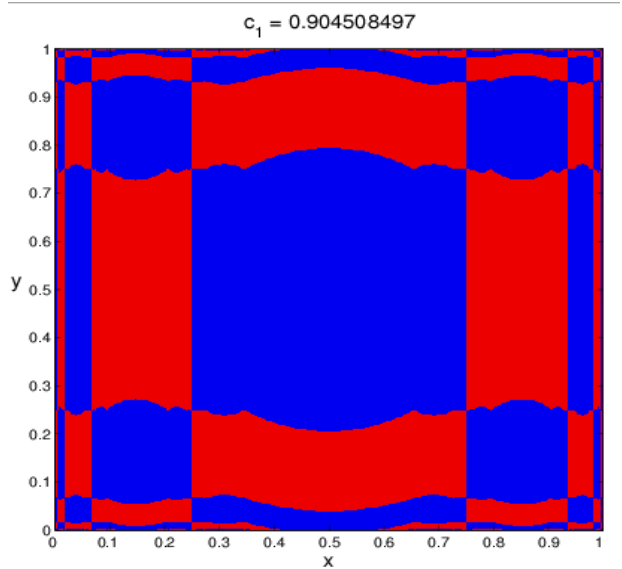,height=3in,width=3.2in}}}
\caption{Basins of attraction for $N=2$ and $c_1 = 0.904508497$. Here $c_1$
is at the upper end of the range for which period 2 attractors 
exist. Note that the boundaries of the basins appear to be non-smooth at some points.}
\label{figc9045}
\end{figure}

\begin{figure}
\centerline{\hbox{\psfig{figure=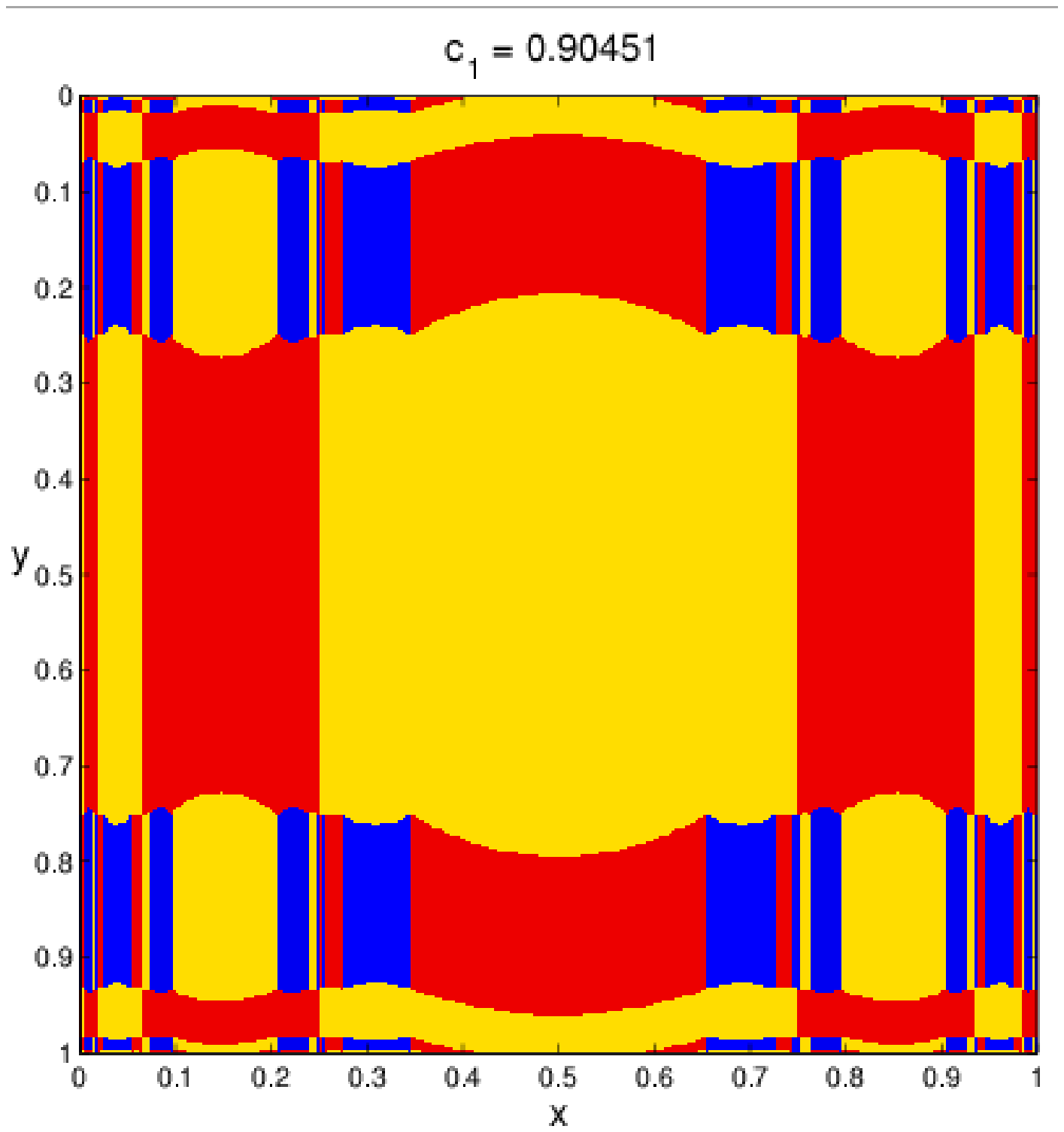,height=3in,width=3.2in}}}
\caption{Basins of attraction for $N=2$ and $c_1 = 0.90451$. This
is just past the range for which period 2 attractors exist. 
Here we observe accumulation of components of the basins in the interior of
$[0,1]^2$. }
\label{figc90451}
\end{figure}

\begin{figure}
\centerline{\hbox{\psfig{figure=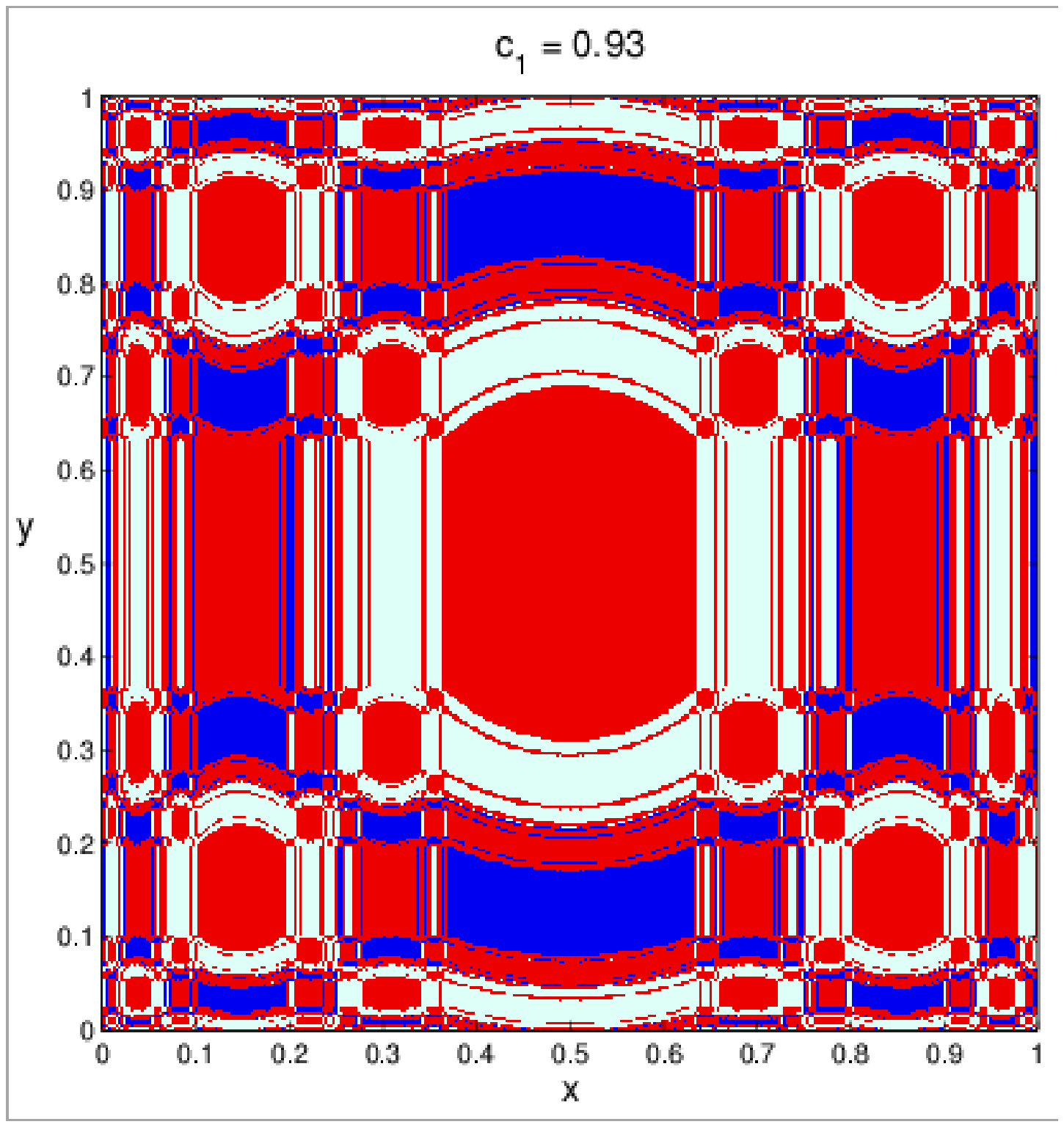,height=3in,width=3.2in}}}
\caption{Basins of attraction for $N=2$ and $c_1 = 0.93$.
Here the accumulation of components of the basins appears to occur 
on a Cantor set even though $c_1 < \xi_2$.}
\label{figc93}
\end{figure}

\begin{figure}
\centerline{\hbox{\psfig{figure=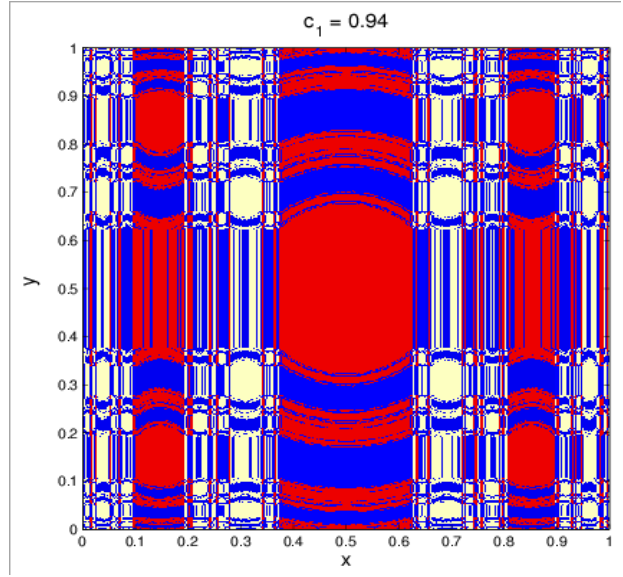,height=3in,width=3.2in}}}
\caption{Basins of attraction for $N=2$ and $c_1 = 0.94 > \xi_2$.}
\label{figc94}
\end{figure}

\begin{figure}
\centerline{\hbox{\psfig{figure=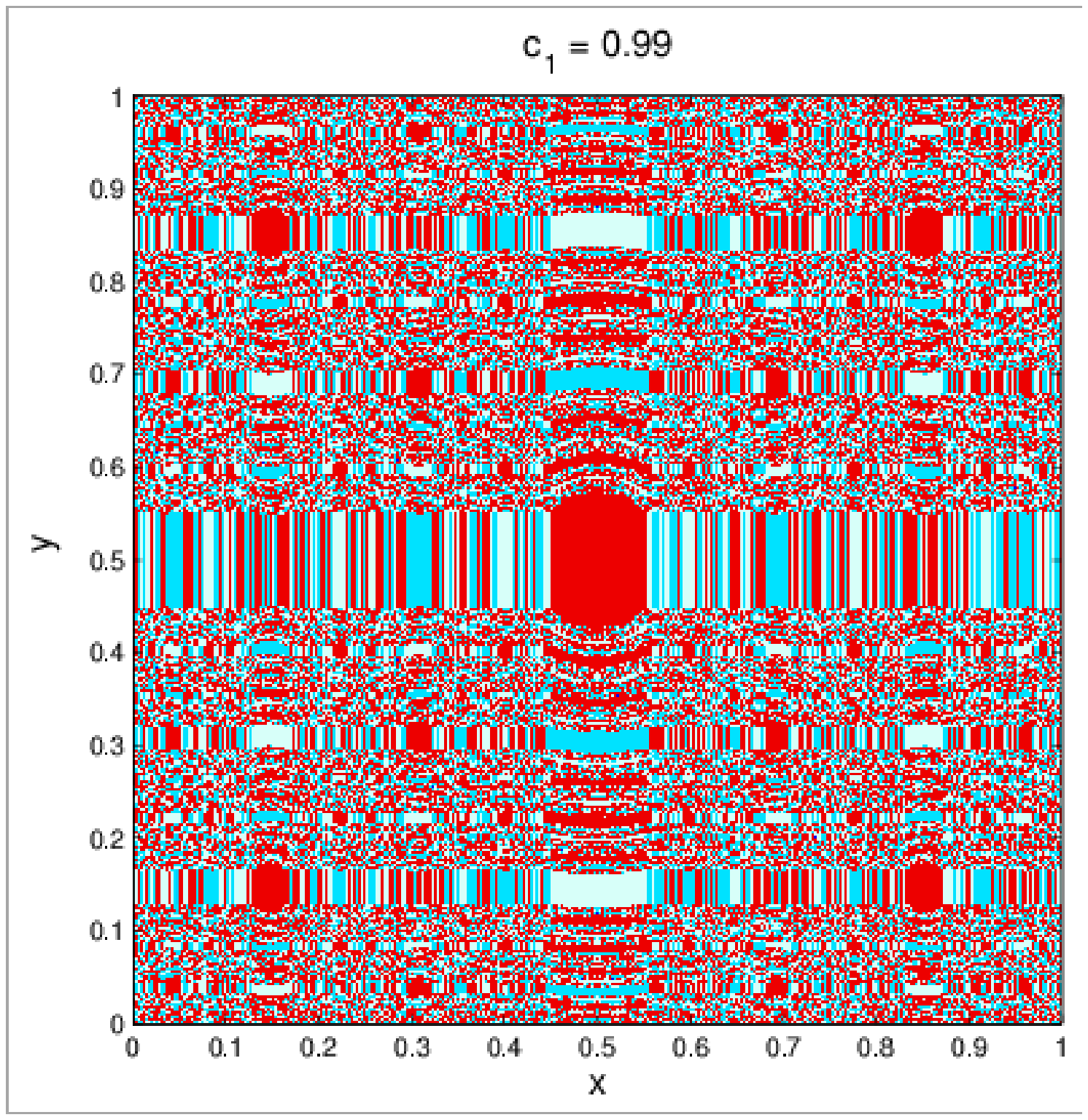,height=3in,width=3.2in}}}
\caption{Basins of attraction for $N=2$ and $c_1 = 0.99$. }
\label{figc99}
\end{figure}

In Figures~\ref{figc84} - \ref{figc99} we show basins of attractors of
various threshold levels for a two site ($N=2$) cascading system.
In these plots we have determined the basins of attractions for each
point in a $499 \times 499$ grid. For each initial condition
the system was allowed to evolve for 100 iterations in order to let
it reach a steady state. Then the sum of the excesses were computed 
for the next 12 iterations. The plots shown are color 
representations of the resulting sum at each point in the grid.
Blue indicates the lowest sum and red indicates the highest.
Different colors must belong to different basins.

Let $C^{-j}$ denote the $j$-th inverse image under $f_{c_1}$ of $C$ and let
$R_j$ denote the compliment of $\bigcup_{i=0}^j C^{-j}$ in $[0,1]$,
i.e. the set of points that avoid $C$ for $j$ steps.
If we let $m_T$ denote the measure induced by the conjugacy
with the tent map, then
$$
    m_T(R_j) = d_1^{j+1}
$$
where $d_1$ is the value of the threshold in the tent map corresponding to
$c_1$. It follows from Theorem 3.1 in \cite{HomYou02} that
$$
      \frac{m(R_j)}{m_T(R_j)} \rightarrow 1
$$
as $j \rightarrow \infty$, where $m$ is Lebesgue measure.
In other words, the asymptotic portion of points falling into
$C$ is the same with respect to Lebesgue measure as with respect
to the induced measure.
Explicitly,  we have that the conjugacy takes $c_1$ to $d_1$ by:
$$
 d_1 = \frac{1}{\pi} \cos^{-1}(1 -2c_1).
$$
From this we may approximate  $m(R_j)$ for various
values of $c_1$ and $j$. For instance when
$c_1 = .9$,we have $d_1 = .7952$ and so after 10 iterations 99\%
of all initial points fall into $C$. For $c_1 = .99$ it
takes 35 iterations for the same percentage of points to 
fall into $C$. After 100 iterations, 99.9\% of initial points
fall into $C$ when $c_1 = .99$. Thus, we expect that 100
iterations in our studies is sufficient for the system to 
reach steady state.

For $3/4< c_1 < .83627...$ in-phase coupling is the only attractor.
For $.83627...<c_1<.90451...$ both in-phase and anti-phase period 2
attractors exists. In Figure~\ref{figc84} we show the basins for $c_1 = 0.84$.
Basins are observed to accumulate at the corners of $[0,1]^2$.

At $c_1 = .904508497...$ a bifurcation occurs; for $c_1$ greater than this value
the period 2 orbit ceases to exist. At the bifurcation 
(Figure~\ref{figc9045}) we notice that the boundary between basin 
is becoming non-smooth at some points. 

In Figure~\ref{figc90451} we find that 
the basins become complicated immediately past this bifurcation. However, from 
the numerics, is not yet fully complex. In contrast the accumulation of components 
of the basins appears to occur on a very complicated set, in 
Figure~\ref{figc93}, for $c_1 = 0.93$ which is before $\xi_2$.

%\begin{figure}
%\centerline{\hbox{\psfig{figure=c95.eps,height=3in,width=3.2in}}}
%\caption{Basins of attraction for $N=2$ and $c_1 = 0.95$. }
%\label{figc95}
%\end{figure}

%\begin{figure}
%\centerline{\hbox{\psfig{figure=c96.eps,height=3in,width=3.2in}}}
%\caption{Basins of attraction for $N=2$ and $c_1 = 0.96$. }
%\label{figc96}
%\end{figure}

%\begin{figure}
%\centerline{\hbox{\psfig{figure=c97.eps,height=3in,width=3.2in}}}
%\caption{Basins of attraction for $N=2$ and $c_1 = 0.97$.}
%\label{figc97}
%\end{figure}

%\begin{figure}
%\centerline{\hbox{\psfig{figure=c98.eps,height=3in,width=3.2in}}}
%\caption{Basins of attraction for $N=2$ and $c_1 = 0.98$.}
%\label{figc98}
%\end{figure}

We find that the basins do not change radically at the bifurcation
value $\xi_2$. In Figures~\ref{figc93}-\ref{figc99} we show basins
for various values of $c_1$. We observe that although the basins are 
quite complicated, they evolve is a very regular way.

Note that points in the basin component in the center of the graph
correspond to in-phase locking. In all our experiments we see that 
this component includes $C^2$, but also extends in the $y$ direction.
This is easily explained by the cascading. For initial points with
$x$ close to $0.5$ and $y$ close to $C$, the cascading pushes 
the next value of $y$ above the threshold. In higher dimensions this
effect will be more pronounced. In fact we have the following:
\begin{prop}
Suppose that for some $j <N$ we have $j (1-c_1) > 1$, then the
central component of in-phase locking will extend to the boundary
of $[0,1]^N$.
\end{prop}

Without cascading, the basins of attraction would be cross-products
of inverse images of $C$, i.e.\ they would be rectangular boxes in $[0,1]$.
The cascading effective perturbs the shape of these basins in the variable
after the first. The further along the array, the greater the potential perturbation.

A summary of our results and numerical observations is contained in Table~\ref{table2}
We have established that cascading maps have multiple basins
for a relatively large range of threshold values.
We have shown that if a cascading system has multiple basins
then those basins have infinitely many components which accumulate at 
the boundary and perhaps also at points in 
the interior. For threshold values above 
$\xi_2 = \frac{2 + \sqrt{3}}{4} \approx 0.933012702$
this  can be very complicated since it can occur at all points
in an Cantor set.

\begin{table}
\centerline{
\begin{tabular}{|c|c|c|c|} \hline
   $c_1$             & Attractors                    & $\Lambda$    & Accumulation  \\ \hline
.75 --  .8362...     & period 2, in-phase             & $\{3/4\}^N$  & N/A \\
.8362 -- .9045...    & in-phase and ripple: $2^{N-1}$  & $\{3/4\}^N$  &  at corners \\
.9045... -- .9330... & various*                         & expanding &  in interior*       \\
.9330... -- .9375    & various*                         & Cantor set & at $\Lambda$      \\ 
.9375... -- 1.0     & in-phase and ripple: $p^{N-1}$  & Cantor set & at $\Lambda$ \\ \hline
\end{tabular}
           }
\label{table2}
\caption{Summary of dynamics for cascading systems. An * indicates
observations from numerics. All others are from proofs.}
\end{table}

The observations suggests that these maps do not
seem promising for applications in pattern recognition, since one
has no hope of approximating an arbitrarily shaped set by one of the basins
of attraction. However, trials using these maps to classify real
medical data have been successful \cite{Par02}. Explanations of this
need to be investigated.

\end{document}